\documentclass[12pt]{article}
\usepackage{amssymb}
\usepackage[all,arc]{xy}
\usepackage{mathrsfs}
\usepackage{amsfonts}
\usepackage{amsmath,amsthm,amsfonts,amssymb,amscd,color}

\allowdisplaybreaks[4]
\newtheorem{defn}{Definition}[section]
\newtheorem{them}[defn]{Theorem}
\newtheorem{lem}[defn]{Lemma}

\newtheorem{cor}[defn]{Corollary}

\newtheorem{con}[defn]{Conjecture}

\numberwithin{equation}{section}
\newenvironment {Proof} {\noindent {\bf Proof.}}{\quad $\square$\par\vspace{3mm}}

\begin{document}
\title{A characterization of arithmetic functions satisfying $f(u^{2}+kv^{2})=f^{2}(u)+kf^{2}(v)$\footnote{L. You's research is supported by the National Natural Science Foundation of China
(Grant No. 11571123) and the Guangdong Provincial Natural Science Foundation (Grant No. 2015A030313377),
Y. Chen's research is supported by the Scientific Research Foundation of Graduate School of South China Normal University (Grant No. 2015lkxm19),
P. Yuan's research is supported by the NSF of China (Grant No. 11271142).}}

\author{ Lihua You\footnote{{\it{Email address:\;}}ylhua@scnu.edu.cn.} 
 \quad Yafei Chen\footnote{{\it{Email address:\;}}764250280@qq.com. } 
 \quad Pingzhi Yuan\footnote{{\it{Corresponding author:\;}}yuanpz@scnu.edu.cn.}
 \quad Jieqin Shen\footnote{{\it{Email address:\;}}648693288@qq.com. }
 }\vskip.2cm
\date{{\small
School of Mathematical Sciences, South China Normal University,\\
Guangzhou, 510631, P.R. China\\
}}
\maketitle
\noindent {\bf Abstract }
In this paper, we mainly discuss the characterization of a class of arithmetic functions $f: N \rightarrow C$ such that $f(u^{2}+kv^2)=f^{2}(u)+kf^{2}(v)$ $(k, u, v \in N)$. We  obtain a characterization with given condition,
propose a conjecture and show the result holds  for $k \in \{2, 3, 4, 5 \}$.

{\it \noindent {\bf Keywords:}}    Arithmetic function; function equation; characterization

{\it \noindent {\bf AMS Classification: } 11A25}

\section{Introduction}
\hskip.6cm Let $C$ denote the set of all complex numbers and $N$ be the set of all  positive integers. An arithmetic function is defined as $f: N \rightarrow C$, which is an important research branch of number theory. A arithmetic function defined on $N$ is called multiplicative if $f(mn)=f(m)f(n)$ for all coprime $m,n \in N$, and is called completely multiplicative if $f(mn)=f(m)f(n)$ for all $m,n \in N$. In recent years, a lot of work on arithmetic function satisfying some Cauchy-like functional equation have been done by researchers.

As was in \cite{1992S}, Spiro showed that if the multiplicative function $f: N \rightarrow C$ satisfies $f(p+q)=f(p)+f(q)$ for all primes $p,q$ and there exists $n_0 \in N$ such that $f(n_0) \neq 0$, then $f(n)=n$ for all $n$. Later, Fang \cite{2011F} extended the conclusion to the equation $f(p+q+r)=f(p)+f(q)+f(r)$. Dubickas et al. \cite{2013D} improved the conclusion to general case $f(p_{1}+p_{2}+ \cdots +p_{k})=f(p_{1})+f(p_{2})+ \cdots +f(p_{k})$, where $k \geq 2$ is fixed, and Chen et al. \cite{2016C} studied the multiplicative function $f$ satisfies $f(p+q+n_{0})=f(p)+f(q)+f(n_{0})$, where $p,q$ are prime and $n_{0}$ is a positive integer. Pong \cite{2006P} considered that if the completely multiplicative function satisfies $f(p+q+pq)=f(p)+f(q)+f(pq)$, $p, q$ are primes, and there exists some $p_{0}$ such that $f(p_{0}) \neq 0$, then $f(n)=n$ for all $n$. De Koninck et al. \cite{1997D} described that if the multiplicative function satisfies $f(1)=1$ and $f(p+m^{2})=f(p)+f(m^{2})$ when $p$ is a prime, then $f(n)=n$ for all $n$. Chung \cite{1996C} characterized the multiplicative and completely multiplicative functions satisfy $f(u^{2}+v^2)=f^{2}(u)+f^{2}(v)$. Phong \cite{2004P} proved that for any $ k \in N$, if the multiplicative function satisfies $f(u^{2}+v^2+k)=f(u^{2})+f(v^2+k)$, $gcd(m, 2k)=1$ and $f(4)f(9) \neq 0$, then $f(n)=n$ for all $n$. Indlekofer et al. \cite{2006I} also discussed the multiplicative function satisfies $f(u^{2}+v^2+k+1)=f(u^{2}+k)+f(v^2+1)$, $f(2) \neq 0$ and $f(5) \neq 1$, then $f(n)=n$ for all $n$, when $gcd(n,2)=1$ holds.

Recently, Ba\v{s}i\'{c} \cite{2014B} considered  a modification of the functional equation treated by Chung: $f(u^{2}+v^2)=f^{2}(u)+f^{2}(v)$, and proved that all such functions can be grouped into three families, namely $f(n) \equiv 0, f(n)=\pm n, f(n)=\pm \frac{1}{k+1}$, and the sign is subjected to the precondition. In this paper, we mainly consider the general cases: $f(u^{2}+kv^2)=f^{2}(u)+kf^{2}(v)$,  we  obtain a characterization with given condition, propose a conjecture and show the result holds  for $k \in \{2, 3, 4, 5 \}$.

\section{Main result and problem}
\hskip.6cm Let $k,u,v\in N$, in order to characterize the function  satisfying $f(u^{2}+kv^{2})=f^{2}(u)+kf^{2}(v)$,
we prove the following theorem firstly.
\begin{them}\label{them201} Let $k\geq2$  be  a positive integer and $A=\left\{\begin{array}{ll}
 6, & k=2; \\
 7, & k=3; \\
 2k, & k\geq 4,
            \end{array}\right.$ $f: N \rightarrow C$ satisfy
\begin{equation}\label{eq201}f(u^{2}+kv^2)=f^{2}(u)+kf^{2}(v)
\end{equation}
for all $u,v\in N$. If one of the following holds for all $n$ $(1\leq n\leq A)$:

\noindent {\rm(1) } $f(n) \equiv 0$;

\noindent{\rm(2) } $f(n) =\left\{\begin{array}{ll}
 n, & \mbox { if  there exist } u, v \in N \mbox { such that } n=u^{2}+kv^{2}; \\
 \pm n, &  \mbox { otherwise},
            \end{array}\right.$

\noindent{\rm(3) } $f(n) =\left\{\begin{array}{ll}
 \frac{1}{k+1}, & \mbox { if  there exist } u, v \in N \mbox { such that } n=u^{2}+kv^{2}; \\
 \pm \frac{1}{k+1}, &  \mbox { otherwise},
            \end{array}\right.$

\noindent Then the corresponding one of {\rm(1) }, {\rm(2) }, {\rm(3) } holds for all $n\in N$.
\end{them}
\begin{Proof}
 Let $a, b, c, d$ be integers with $ab \pm kcd > 0$ and $ad \pm bc > 0$. Firstly, we have
\begin{equation}\label{eq202}
(ab+kcd)^{2}+k(ad-bc)^{2}=(ab-kcd)^{2}+k(ad+bc)^{2}.\end{equation}
Applying $f$ to both sides of (\ref{eq202}) and using equation (\ref{eq201}) we have
\begin{equation}\label{eq203} f^{2}(ab+kcd)+kf^{2}(ad-bc)=f^{2}(ab-kcd)+kf^{2}(ad+bc).
\end{equation}

Now we only need to show one of (1), (2), (3) holds when $n>A$. By induction on $n$, we complete the proof by the following two cases.

\noindent {\bf Case 1: } $k$ is odd.

Then $k \geq 3$.

\noindent {\bf Subcase 1.1:} $n$ is odd.

Let $n=2l+1$, then  by $n > A$ we have  $l\geq \left\{ \begin{array}{cc}
                                                 4, & \mbox  { if } k=3; \\
                                                 k, &  \mbox  { if } k\geq5.
                                               \end{array}\right.$
   We take $a=l-\frac{k-1}{2}, b=2, c=d=1$ in  (\ref{eq203}), then $2l-2k+1>0,$ $l-\frac{k-1}{2}-2>0$ and
\begin{equation}\label{eq204} f^{2}(n)=f^{2}(2l-2k+1)+kf^{2}(l-\frac{k-1}{2}+2)-kf^{2}(l-\frac{k-1}{2}-2).\end{equation}

\noindent \textbf{Subcase 1.1.1: } {\rm(1) } holds for $ 1\leq n\leq A$.

By the inductive assumption and $f(n)\equiv 0$ for  $ 1\leq n\leq A $, we have
$f^{2}(2l-2k+1)=f^{2}(l-\frac{k-1}{2}+2)=f^{2}(l-\frac{k-1}{2}-2)=0$.
Then   $f^{2}(n)=0$ by (\ref{eq204}) and thus $f(n)=0$.

\noindent \textbf{Subcase 1.1.2: } {\rm(2) } holds for $ 1\leq n\leq A$.

If  there exist $u, v \in N$  such that $n=u^{2}+kv^{2}$, then
by the inductive assumption and {\rm(2) } holds for $ 1\leq n\leq A$,
we have $f^{2}(u)=f^{2}(v)=\frac{1}{(k+1)^{2}}$, and thus $f(n)=\frac{1}{k+1}$ by (\ref{eq201}).

Otherwise, by the inductive assumption and {\rm(2) } holds for $ 1\leq n\leq A$, we have $f^{2}(2l-2k+1)=f^{2}(l-\frac{k-1}{2}+2)=f^{2}(l-\frac{k-1}{2}-2)=\frac{1}{(k+1)^{2}}$. Then $f(n)=\pm\frac{1}{k+1}$ by (\ref{eq204}).

\noindent \textbf{Subcase 1.1.3: } {\rm(3) }holds for $ 1\leq n\leq A$.

If there exist $u, v \in N$  such that $n=u^{2}+kv^{2}$, then
by the inductive assumption, we have $f^{2}(u)=u^2$, $f^{2}(v)=v^2$, thus $f(n)=n$ by (\ref{eq201}).

Otherwise, by the inductive assumption,
we have  $f^{2}(2l-2k+1)=(2l-2k+1)^{2}$, $f^{2}(l-\frac{k-1}{2}+2)=(l-\frac{k-1}{2}+2)^{2}$ and $f^{2}(l-\frac{k-1}{2}-2)=(l-\frac{k-1}{2}-2)^{2}$. Then we have $f^{2}(n)=(2l-2k+1)^{2}+k(l-\frac{k-1}{2}+2)^{2}-k(l-\frac{k-1}{2}-2)^{2}=4l^{2}+4l+1=(2l+1)^{2}=n^2$ by (\ref{eq204}), and thus $f(n)=\pm n$.

Combining the above three subcases, we complete the proof when $n$ is odd.

\noindent {\bf Subcase 1.2:} $n$ is even.

Let $n=2l$, then $l \geq k+1$ by $n > A$. We take $a=2l-k, b=c=d=1$ in  (\ref{eq203}), then $2l-2k>0$, $2l-k-1>0$ and
\begin{equation}\label{eq205} f^{2}(n)=f^{2}(2l-2k)+kf^{2}(2l-k+1)-kf^{2}(2l-k-1).\end{equation}

\noindent \textbf{Subcase 1.2.1: } {\rm(1) } holds for $ 1\leq n\leq A$.

By the inductive assumption and $f(n)\equiv 0$ for  $ 1\leq n\leq A $, we have
$f^{2}(2l-2k)=f^{2}(2l-k+1)=f^{2}(2l-k-1)=0$.
Then   $f^{2}(n)=0$ by (\ref{eq205}) and thus $f(n)=0$.

\noindent \textbf{Subcase 1.2.2: } {\rm(2) } holds for $ 1\leq n\leq A$.

If  there exist $u, v \in N$  such that $n=u^{2}+kv^{2}$, then
by the inductive assumption and {\rm(2) } holds for $ 1\leq n\leq A$,
we have $f^{2}(u)=f^{2}(v)=\frac{1}{(k+1)^{2}}$, and thus $f(n)=\frac{1}{k+1}$ by (\ref{eq201}).

Otherwise, by the inductive assumption and {\rm(2) } holds for $ 1\leq n\leq A$, we have $f^{2}(2l-2k)=f^{2}(2l-k+1)=f^{2}(2l-k-1)=\frac{1}{(k+1)^{2}}$. Then $f(n)=\pm\frac{1}{k+1}$ by (\ref{eq205}).

\noindent \textbf{Subcase 1.2.3: } {\rm(3) }holds for $ 1\leq n\leq A$.

If there exist $u, v \in N$  such that $n=u^{2}+kv^{2}$, then
by the inductive assumption, we have $f^{2}(u)=u^2$, $f^{2}(v)=v^2$, thus $f(n)=n$ by (\ref{eq201}).

Otherwise, by the inductive assumption,
we have  $f^{2}(2l-2k)=(2l-2k)^{2}$, $f^{2}(2l-k+1)=(2l-k+1)^{2}$ and $f^{2}(2l-k-1)=(2l-k-1)^{2}$. Then we have $f^{2}(n)=(2l-2k)^{2}+k(2l-k+1)^{2}-k(2l-k-1)^{2}=4l^{2}=n^2$ by (\ref{eq205}), and thus $f(n)=\pm n$.

Combining the above three subcases, we complete the proof when $n$ is even.

\noindent {\bf Case 2: } $k$ is even.

\noindent {\bf Subcase 2.1:} $n$ is odd.

Let $n=2l+1$, then  by $n > A$ we have  $l\geq \left\{ \begin{array}{cc}
                                                 3, & \mbox  { if } k=2; \\
                                                 k, &  \mbox  { if } k\geq4.
                                               \end{array}\right.$
 We take $a=2l-k+1, b=c=d=1$ in equation (\ref{eq203}), then $2l-2k+1>0$, $2l-k>0$ and
\begin{equation}\label{eq206} f^{2}(n)=f^{2}(2l-2k+1)+kf^{2}(2l-k+2)-kf^{2}(2l-k). \end{equation}

\noindent {\bf Subcase 2.2:} $n$ is even.

Let $n=2l$,  then  by $n > A$ we have  $l\geq \left\{ \begin{array}{cc}
                                                 4, & \mbox  { if } k=2; \\
                                                 k+1, &  \mbox  { if } k\geq4.
                                               \end{array}\right.$
 We take $a=l-\frac{k}{2}, b=c=d=1$ in equation (\ref{eq203}), then $2l-2k>0$, $l-\frac{k}{2}-2>0$ and

\begin{equation}\label{eq207} f^{2}(n)=f^{2}(2l-2k)+kf^{2}(l-\frac{k}{2}+2)-kf^{2}(l-\frac{k}{2}-2). \end{equation}

Similar to the proof of Case 1, we can complete the proof by the  assumption, (\ref{eq201}), (\ref{eq206}) or (\ref{eq207}), so we omit it.
\end{Proof}

Based on  the result of Theorem \ref{them201}, we  propose the following conjecture for further research.
\begin{con}\label{con202} Let $k\in N$, $f: N \rightarrow C$ satisfy (\ref{eq201})
for all $u,v\in N$. Then one of the following holds:

\noindent {\rm(1) } $f(n) \equiv 0$;

\noindent{\rm(2) } $f(n) =\left\{\begin{array}{ll}
 n, & \mbox { if  there exist } u, v \in N \mbox { such that } n=u^{2}+kv^{2}; \\
 \pm n, &  \mbox { otherwise}.
            \end{array}\right.$

\noindent{\rm(3) } $f(n) =\left\{\begin{array}{ll}
 \frac{1}{k+1}, & \mbox { if  there exist } u, v \in N \mbox { such that } n=u^{2}+kv^{2}; \\
 \pm \frac{1}{k+1}, &  \mbox { otherwise}.
            \end{array}\right.$
\end{con}

In the following sections, we will show Conjecture \ref{con202} holds for the cases $k=2,3,4,5$.

\section{The proof of $k=2$}
\hskip.6cm In this section, we will prove Conjecture \ref{con202} holds for $k=2$.

\begin{lem}\label{lem302}Let $f: N \rightarrow C$ satisfy
\begin{equation}\label{eq301}f(u^{2}+2v^2)=f^{2}(u)+2f^{2}(v)
\end{equation}
for all $u,v\in N$,  $f(1)=a$ and $f(2)=b$.
Then we have
\begin{equation}\label{eq302}f^{2}(3)=9a^{4},
\end{equation}
\begin{equation}\label{eq303}f^{2}(4)=\frac{1}{2}(27a^{4}-3a^{2}+2b^{2}),
\end{equation}
\begin{equation}\label{eq304}f^{2}(5)=27a^{4}-2a^{2},
\end{equation}
\begin{equation}\label{eq305}f^{2}(6)=36a^{4}-4a^{2}+b^{2}.
\end{equation}
\end{lem}

\begin{Proof}
Since $f(1)=a, f(2)=b$ and $3=1^{2}+2\times1^{2}$, then we have $f(3)=3a^{2}$ and thus (\ref{eq302}) holds.
Noting that $27=3^{2}+2\times3^{2}=5^{2}+2\times1^{2}, 33=1^{2}+2\times4^{2}=5^{2}+2\times2^{2}$ and $54=2^{2}+2\times5^{2}=6^{2}+2\times3^{2}$, then we have
\begin{equation}\label{eq306}
\left\{\begin{array}{ll}
f^{2}(3)+2f^{2}(3)=f^{2}(5)+2f^{2}(1), \\
f^{2}(1)+2f^{2}(4)=f^{2}(5)+2f^{2}(2), \\
f^{2}(2)+2f^{2}(5)=f^{2}(6)+2f^{2}(3).
\end{array}\right.
\end{equation}
Solving (\ref{eq306}), we have (\ref{eq303})-(\ref{eq305}) hold, which complete the proof.
\end{Proof}

Now we evaluate $f(1)$.
\begin{lem}\label{lem303} $f(1)\in\{0,1,-1, \frac{1}{3}, -\frac{1}{3}\}$.
\end{lem}
\begin{Proof}
Firstly, by $f(1)=a, f(2)=b$ and $f(6)=f(2^{2}+2\times1^{2})=f^{2}(2)+2f^{2}(1)=b^{2}+2a^{2}$, we have
\begin{equation}\label{eq307} f^{2}(6)=(b^{2}+2a^{2})^{2}.
\end{equation}

Similarly, by $f(9)=f(1^{2}+2\times2^{2})=f^{2}(1)+2f^{2}(2)=a^{2}+2b^{2}$, we have
\begin{equation}\label{eq308} f^{2}(9)=(a^{2}+2b^{2})^{2}.
\end{equation}

Noting that $51=1^{2}+2\times5^{2}=7^{2}+2\times1^{2}$ and $99=1^{2}+2\times7^{2}=9^{2}+2\times3^{2}$, then we have
 $f^{2}(1)+2f^{2}(5)=f^{2}(7)+2f^{2}(1),$ $f^{2}(1)+2f^{2}(7)=f^{2}(9)+2f^{2}(3)$ and thus
 \begin{equation}\label{eq309} f^{2}(9)=90a^{4}-9a^{2}.
\end{equation}

Combining (\ref{eq305}), (\ref{eq307}), (\ref{eq308}) and (\ref{eq309}), we have
\begin{equation*}
\left\{\begin{array}{ll}
  (b^{2}+2a^{2})^{2}=36a^{4}-4a^{2}+b^{2}, \\
  (a^{2}+2b^{2})^{2}=90a^{4}-9a^{2}.
\end{array}\right.
\end{equation*}
Then
\begin{equation}\label{eq310}
b^{2}=\frac{39a^{4}-7a^{2}}{12a^{2}-4},
\end{equation}
and  $(\frac{39a^{4}-7a^{2}}{12a^{2}-4})^{2}+\frac{4a^{2}(39a^{4}-7a^{2})}{12a^{2}-4}+4a^{4}=36a^{4}-4a^{2}+\frac{39a^{4}-7a^{2}}{12a^{2}-4}$.
By simplifying this equation, we have
$a^{2}(a-1)(a+1)(3a-1)(3a+1)(15a^{2}-4)=0.$
Thus $$a_{1}=a_{2}=0, a_{3}=1, a_{4}=-1, a_{5}=\frac{1}{3}, a_{6}=-\frac{1}{3}, a_{7}=\frac{2\sqrt{15}}{15}, a_{8}=-\frac{2\sqrt{15}}{15}.$$

Now we show that $f(1)=\pm \frac{2\sqrt{15}}{15}$ is impossible. Otherwise, we have $a^{2}=\frac{4}{15}$
and $b^{2}=-\frac{17}{15}$ by  (\ref{eq310}).
Since $f(12)=f(2^{2}+2\times2^{2})=f^{2}(2)+2f^{2}(2)=3b^{2}$, then we have
\begin{equation}\label{eq311} f^{2}(12)=9b^{4}.
\end{equation}

On the other hand, we note  that $108=10^{2}+2\times2^{2}=6^{2}+2\times6^{2}$ and $216=4^{2}+2\times10^{2}=12^{2}+2\times6^{2}$, then we have
 $ f^{2}(10)+2f^{2}(2)=f^{2}(6)+2f^{2}(6),$ $f^{2}(4)+2f^{2}(10)=f^{2}(12)+2f^{2}(6)$ and
\begin{equation}\label{eq312} f^{2}(12)=f^{2}(4)+4f^{2}(6)-4b^{2}=\frac{315a^{4}-35a^{2}+2b^{2}}{2}.
\end{equation}

By (\ref{eq311})-(\ref{eq312}), we have
$$\frac{315a^{4}-35a^{2}+2b^{2}}{2}-9b^{4}=0.$$
Unfortunately,  when we put $a^{2}=\frac{4}{15}$ and $b^{2}=-\frac{17}{15}$ in $\frac{315a^{4}-35a^{2}+2b^{2}}{2}-9b^{4}$,
we obtain $\frac{315a^{4}-35a^{2}+2b^{2}}{2}-9b^{4}\not=0$. Thus it   implies a contradiction and we complete the proof.
\end{Proof}

By Lemmas \ref{lem302}-\ref{lem303} and direct calculation, we obtain the following Coroll1ary.
\begin{cor}\label{cor304} Let $1\leq n \leq6$ and $f$ satisfy the function equation (\ref{eq301}). Then one of the following holds.

\noindent {\rm(1) } $f(n) \equiv 0$;

\noindent {\rm(2) } $f(n) =\left\{\begin{array}{ll}
 n, & \mbox { if } n\in \{3,6\}; \\
 \pm n, &  \mbox { otherwise}.
            \end{array}\right.$

\noindent {\rm(3) }$f(n) =\left\{\begin{array}{ll}
 \frac{1}{3}, & \mbox { if } n\in \{3,6\}; \\
 \pm \frac{1}{3}, &  \mbox { otherwise}.
            \end{array}\right.$
\end{cor}

Now, by Theorem \ref{them201} and Corollary \ref{cor304}, we obtain Theorem \ref{them301} directly.

\begin{them}\label{them301}
Conjecture \ref{con202} holds for $k=2$.
\end{them}

\section{The proof of $k=3$}
\hskip.6cm In this section, we will prove Conjecture \ref{con202} holds for $k=3$.

\begin{lem}\label{lem402}Let $f: N \rightarrow C$ satisfy
\begin{equation}\label{eq401}f(u^{2}+3v^2)=f^{2}(u)+3f^{2}(v)
\end{equation}
for all $u,v\in N$,  $f(1)=a$ and $f(2)=b$.
Then we have
\begin{equation}\label{eq402}f^{2}(3)=\frac{1}{3}(8b^{2}-5a^{2}), \end{equation}
\begin{equation}\label{eq403}f^{2}(4)=5b^{2}-4a^{2}, \end{equation}
\begin{equation}\label{eq404}f^{2}(5)=8b^{2}-7a^{2}, \end{equation}
\begin{equation}\label{eq405}f^{2}(6)=\frac{1}{3}(35b^{2}-32a^{2}), \end{equation}
\begin{equation}\label{eq406}f^{2}(7)=16b^{2}-15a^{2}, \end{equation}
\end{lem}

\begin{Proof}
Since $f(1)=a, f(2)=b$ and $28=5^{2}+3\times1^{2}=1^{2}+3\times3^{2}=4^{2}+3\times2^{2}$ and $52=2^{2}+3\times4^{2}=5^{2}+3\times3^{2}$,
we have
\begin{equation}\label{eq407}
\left\{\begin{array}{ll}
f^{2}(5)=3f^{2}(3)-2a^{2}, \\
f^{2}(4)=3f^{2}(3)+a^{2}-3b^{2},\\
b^{2}+3f^{2}(4)=f^{2}(5)+3f^{2}(3).
\end{array}\right.
\end{equation}
Solving (\ref{eq407}), we have (\ref{eq402}), (\ref{eq403}) and (\ref{eq404}) hold.

Similarly, it is easy to find that $84=6^{2}+3\times4^{2}=3^{2}+3\times5^{2}$ and $124=7^{2}+3\times5^{2}=4^{2}+3\times6^{2}$, then we have
\begin{equation}\label{eq408}
\left\{\begin{array}{ll}
f^{2}(6)+3f^{2}(4)=f^{2}(3)+3f^{2}(5), \\
f^{2}(7)+3f^{2}(5)=f^{2}(4)+3f^{2}(6).
\end{array}\right.
\end{equation}
Solving (\ref{eq408}), we have (\ref{eq405}) and (\ref{eq406}) hold. Then we complete the proof.
\end{Proof}

Now we evaluate $f(1)$.
\begin{lem}\label{lem403} $f(1)\in\{0,1,-1, \frac{1}{4}, -\frac{1}{4}\}$.
\end{lem}
\begin{Proof}
Firstly, by $f(1)=a$ and $f(4)=f(1^{2}+3\times1^{2})=f^{2}(1)+3f^{2}(1)=4f^{2}(1)=4a^{2}$, we have
\begin{equation}\label{eq409} f^{2}(4)=16a^{4}.
\end{equation}

Noting that $112=10^{2}+3\times2^{2}=2^{2}+3\times6^{2}$, $208=14^{2}+3\times2^{2}=10^{2}+3\times6^{2}$ and $304=16^{2}+3\times4^{2}=14^{2}+3\times6^{2}$, then we have
\begin{equation*}
\left\{\begin{array}{ll}
  f^{2}(10)+3f^{2}(2)=f^{2}(2)+3f^{2}(6), \\
  f^{2}(14)+3f^{2}(2)=f^{2}(10)+3f^{2}(6), \\
  f^{2}(16)+3f^{2}(4)=f^{2}(14)+3f^{2}(6),
\end{array}\right.
\end{equation*}
and \begin{equation}\label{eq410} f^{2}(16)=9f^{2}(6)-5f^{2}(2)-3f^{2}(4)=85b^{2}-84a^{2}.
\end{equation}

On the other hand, by $f(2)=b$ and $f(16)=f(2^{2}+3\times2^{2})=f^{2}(2)+3f^{2}(2)=4f^{2}(2)=4b^{2}$, we have
\begin{equation}\label{eq411} f^{2}(16)=16b^{4}.
\end{equation}

Combining (\ref{eq403}), (\ref{eq409}), (\ref{eq410}) and (\ref{eq411}), we have
\begin{equation*}
\left\{\begin{array}{ll}
  5b^{2}-4a^{2}=16a^{4}, \\
  85b^{2}-84a^{2}=16b^{4}.
\end{array}\right.
\end{equation*}
Then
\begin{equation}\label{eq412}
b^{2}=\frac{1}{5}(16a^{4}+4a^{2}),
\end{equation}
and  $17(16a^{4}+4a^{2})-84a^{2}=\frac{16}{25}(16a^{4}+4a^{2})^{2}$.
By simplifying this equation, we have
$a^{2}(a-1)(a+1)(4a-1)(4a+1)(16a^{2}+25)=0.$
Thus $$a_{1}=a_{2}=0, a_{3}=1, a_{4}=-1, a_{5}=\frac{1}{4}, a_{6}=-\frac{1}{4}, a_{7}=\frac{5}{4}i, a_{8}=-\frac{5}{4}i.$$

Now we show that $f(1)=\pm \frac{5}{4}i$ is impossible. Otherwise, we have $a^{2}=-\frac{25}{16}$
and $b^{2}=\frac{105}{16}$ by  (\ref{eq412}).
Noting that $f(7)=f(2^{2}+3\times1^{2})=f^{2}(2)+3f^{2}(1)=b^{2}+3a^{2}$, by  (\ref{eq406}), we have
$$(b^{2}+3a^{2})^{2}-(16b^{2}-15a^{2})=0.$$
But when we put $a^{2}=-\frac{25}{16}$ and $b^{2}=\frac{105}{16}$ in $(b^{2}+3a^{2})^{2}-(16b^{2}-15a^{2})$,
we obtain $(b^{2}+3a^{2})^{2}-(16b^{2}-15a^{2})\not=0$ which implies a contradiction and we complete the proof.
\end{Proof}

By Lemmas \ref{lem402}-\ref{lem403} and direct calculation, we obtain the following Coroll1ary.
\begin{cor}\label{cor404} Let $1\leq n \leq7$ and $f$ satisfy the function equation (\ref{eq401}). Then one of the following holds.

\noindent {\rm(1) } $f(n) \equiv 0$;

\noindent {\rm(2) } $f(n) =\left\{\begin{array}{ll}
 n, & \mbox { if } n\in\{4,7\}; \\
 \pm n, &  \mbox { otherwise}.
            \end{array}\right.$

\noindent {\rm(3) }$f(n) =\left\{\begin{array}{ll}
 \frac{1}{4}, & \mbox { if } n\in\{4,7\}; \\
 \pm \frac{1}{4}, &  \mbox { otherwise}.
            \end{array}\right.$
\end{cor}

Now, by Theorem \ref{them201} and Corollary \ref{cor404}, we obtain Theorem \ref{them401} directly.

\begin{them}\label{them401}
Conjecture \ref{con202} holds for $k=3$.
\end{them}

\section{The proof of $k=4$}
\hskip.6cm In this section, we will prove Conjecture \ref{con202} holds for $k=4$.

\begin{lem}\label{lem502}Let $f: N \rightarrow C$ satisfy
\begin{equation}\label{eq501}f(u^{2}+4v^2)=f^{2}(u)+4f^{2}(v)
\end{equation}
for all $u,v\in N$,  $f(1)=a$ and $f(2)=b$.
Then (\ref{eq402})-(\ref{eq406}) and (\ref{eq502}) hold:
\begin{equation}\label{eq502}f^{2}(8)=21b^{2}-20a^{2}. \end{equation}
\end{lem}

\begin{Proof}
Since $f(1)=a, f(2)=b$,  $20=4^{2}+4\times1^{2}=2^{2}+4\times2^{2}$,  $65=7^{2}+4\times2^{2}=1^{2}+4\times4^{2}$  and $68=8^{2}+4\times1^{2}=2^{2}+4\times4^{2}$, we have
   $$f^{2}(4)=5f^{2}(2)-4f^{2}(1)=5b^{2}-4a^{2},$$
  $$f^{2}(7)=f^{2}(1)+4f^{2}(4)-4f^{2}(2)=16b^{2}-15a^{2},$$
  \noindent and $$ f^{2}(8)=f^{2}(2)+4f^{2}(4)-4f^{2}(1)=21b^{2}-20a^{2}.$$

Noting that $200=10^{2}+4\times5^{2}=2^{2}+4\times7^{2}$ and $104=10^{2}+4\times1^{2}=2^{2}+4\times5^{2}$, we have
\begin{equation}\label{eq503}
\left\{\begin{array}{ll}
f^{2}(10)+4f^{2}(5)=f^{2}(2)+4f^{2}(7), \\
f^{2}(10)+4f^{2}(1)=f^{2}(2)+4f^{2}(5).
\end{array}\right.
\end{equation}

Solving (\ref{eq503}), we obtain $f^{2}(10)=33b^{2}-32a^{2}$  and (\ref{eq404}) holds.

Similarly, it is easy  to find that $100=6^{2}+4\times4^{2}=8^{2}+4\times3^{2}$, $265=3^{2}+4\times8^{2}=11^{2}+4\times6^{2}$ and $125=11^{2}+4\times1^{2}=5^{2}+4\times5^{2}$. Then we have
\begin{equation}\label{eq504}
\left\{\begin{array}{ll}
 f^{2}(6)+4f^{2}(4)=f^{2}(8)+4f^{2}(3), \\
 f^{2}(3)+4f^{2}(8)=f^{2}(11)+4f^{2}(6), \\
 f^{2}(11)+4f^{2}(1)=f^{2}(5)+4f^{2}(5).
\end{array}\right.
\end{equation}

Solving (\ref{eq504}), (\ref{eq402}) and (\ref{eq405}) hold, and we complete the proof.
\end{Proof}

Now we evaluate $f(1)$.
\begin{lem}\label{lem503} $f(1)\in\{0,1,-1, \frac{1}{5}, -\frac{1}{5}\}$.
\end{lem}
\begin{Proof}
Firstly, by $f(1)=a$ and $f(5)=f(1^{2}+4\times1^{2})=f^{2}(1)+4f^{2}(1)=5f^{2}(1)=5a^{2}$, we have
\begin{equation}\label{eq505} f^{2}(5)=25a^{4}.
\end{equation}

Similarly, by $f(2)=b$ and $f(20)=f(2^{2}+4\times2^{2})=f^{2}(2)+4f^{2}(2)=5f^{2}(2)=5b^{2}$, we have
\begin{equation}\label{eq506} f^{2}(20)=25b^{4}.
\end{equation}

Noting that $404=20^{2}+4\times1^{2}=2^{2}+4\times10^{2}$, by (\ref{eq503}) we have
\begin{equation}\label{eq507} f^{2}(20)=f^{2}(2)+4f^{2}(10)-4f^{2}(1)=133b^{2}-132a^{2}.
\end{equation}

Combining (\ref{eq404}), (\ref{eq505}), (\ref{eq506}) and (\ref{eq507}), we have
\begin{equation*}
\left\{\begin{array}{ll}
  8b^{2}-7a^{2}=25a^{4}, \\
  133b^{2}-132a^{2}=25b^{4}.
\end{array}\right.
\end{equation*}
Then we have
\begin{equation}\label{eq508}
b^{2}=\frac{1}{8}(25a^{4}+7a^{2}),
\end{equation}
and  $\frac{133}{8}(25a^{4}+7a^{2})-132a^{2}=\frac{25}{64}(25a^{4}+7a^{2})^{2}$.
By simplifying this equation, we have
$a^{2}(a-1)(a+1)(5a-1)(5a+1)(5a^{2}+8)=0.$
Thus $$a_{1}=a_{2}=0, a_{3}=1, a_{4}=-1, a_{5}=\frac{1}{5}, a_{6}=-\frac{1}{5}, a_{7}=\frac{2\sqrt10}{5}i, a_{8}=-\frac{2\sqrt10}{5}i.$$

Now we show that $f(1)=\pm \frac{2\sqrt10}{5}i$ is impossible. Otherwise, we have $a^{2}=-\frac{8}{5}$
and $b^{2}=\frac{33}{5}$ by  (\ref{eq508}).
Noting that $f(8)=f(2^{2}+4\times1^{2})=f^{2}(2)+4f^{2}(1)=b^{2}+4a^{2}$, by  (\ref{eq502}), we have
$$(b^{2}+4a^{2})^{2}-(21b^{2}-20a^{2})=0.$$
Unfortunately,  when we put $a^{2}=-\frac{8}{5}$ and $b^{2}=\frac{33}{5}$ in $(b^{2}+4a^{2})^{2}-(21b^{2}-20a^{2})$,
we obtain $(b^{2}+4a^{2})^{2}-(21b^{2}-20a^{2})\not=0$ which implies a contradiction and we complete the proof.
\end{Proof}

By Lemmas \ref{lem502}-\ref{lem503} and direct calculation, we obtain the following Coroll1ary.
\begin{cor}\label{cor504} Let $1\leq n \leq8$ and $f$ satisfy the function equation (\ref{eq501}). Then one of the following holds.

\noindent {\rm(1) } $f(n) \equiv 0$;

\noindent {\rm(2) } $f(n) =\left\{\begin{array}{ll}
 n, & \mbox { if } n \in \{5,8\}; \\
 \pm n, &  \mbox { otherwise}.
            \end{array}\right.$

\noindent {\rm(3) }$f(n) =\left\{\begin{array}{ll}
 \frac{1}{5}, & \mbox { if } n \in \{5,8\}; \\
 \pm \frac{1}{5}, &  \mbox { otherwise}.
            \end{array}\right.$
\end{cor}

Now, by Theorem \ref{them201} and Corollary \ref{cor504}, we obtain Theorem \ref{them501} directly.

\begin{them}\label{them501}
Conjecture \ref{con202} holds for $k=4$.
\end{them}

\section{The proof of $k=5$}
\hskip.6cm In this section, we will prove Conjecture \ref{con202} holds for $k=5$.

\begin{lem}\label{lem602}Let $f: N \rightarrow C$ satisfy
\begin{equation}\label{eq601}f(u^{2}+5v^2)=f^{2}(u)+5f^{2}(v)
\end{equation}
for all $u,v\in N$,  $f(1)=a$ and $f(2)=b$.
Then (\ref{eq402})-(\ref{eq406}), (\ref{eq502}), (\ref{eq602}) and (\ref{eq603}) hold:
\begin{equation}\label{eq602}f^{2}(9)=\frac{1}{3}(80b^{2}-77a^{2}), \end{equation}
\begin{equation}\label{eq603}f^{2}(10)=33b^{2}-32a^{2}.\end{equation}
\end{lem}

\begin{Proof}
Since $f(1)=a, f(2)=b$ and $21=4^{2}+5\times1^{2}=1^{2}+5\times2^{2}$, $84=8^{2}+5\times2^{2}=2^{2}+5\times4^{2}$,
$69=7^{2}+5\times2^{2}=8^{2}+5\times1^{2}$,  $129=7^{2}+5\times4^{2}=2^{2}+5\times5^{2}$,
$345=10^{2}+5\times7^{2}=5^{2}+5\times8^{2}$, $54=3^{2}+5\times3^{2}=7^{2}+5\times1^{2}$,
$81=6^{2}+5\times3^{2}=1^{2}+5\times4^{2}$, $126=9^{2}+5\times3^{2}=1^{2}+5\times5^{2}$,
we have
\begin{equation}\label{eq604}
\left\{\begin{array}{c}
               f^{2}(4)+5f^{2}(1)=f^{2}(1)+5f^{2}(2), \\
              f^{2}(8)+5f^{2}(2)=f^{2}(2)+5f^{2}(4), \\
               f^{2}(7)+5f^{2}(2)=f^{2}(8)+5f^{2}(1), \\
               f^{2}(7)+5f^{2}(4)=f^{2}(2)+5f^{2}(5), \\
              f^{2}(10)+5f^{2}(7)=f^{2}(5)+5f^{2}(8), \\
              f^{2}(3)+5f^{2}(3)=f^{2}(7)+5f^{2}(1), \\
               f^{2}(6)+5f^{2}(3)=f^{2}(1)+5f^{2}(4), \\
               f^{2}(9)+5f^{2}(3)=f^{2}(1)+5f^{2}(5).
             \end{array}\right.
\end{equation}
Solving (\ref{eq604}),  (\ref{eq402})-(\ref{eq406}), (\ref{eq502}), (\ref{eq602})-(\ref{eq603}) hold, so we complete the proof.
\end{Proof}

Now we evaluate $f(1)$.
\begin{lem}\label{lem603}
$f(1)\in\{0,1,-1, \frac{1}{6}, -\frac{1}{6}\}$.
\end{lem}
\begin{Proof}
Firstly, by $f(1)=a$ and $f(6)=f(1^{2}+5\times1^{2})=f^{2}(1)+5f^{2}(1)=6f^{2}(1)=6a^{2}$, we have
\begin{equation}\label{eq605} f^{2}(6)=36a^{4}.
\end{equation}

Similarly, by $f(2)=b$ and $f(9)=f(2^{2}+5\times1^{2})=f^{2}(2)+5f^{2}(1)=b^{2}+5a^{2}$, we have
\begin{equation}\label{eq606} f^{2}(9)=(b^{2}+5a^{2})^{2}.
\end{equation}

Combining (\ref{eq405}), (\ref{eq602}), (\ref{eq605}) and (\ref{eq606}), we have
\begin{equation*}
\left\{\begin{array}{ll}
  \frac{1}{3}(35b^{2}-32a^{2})=36a^{4} \\
 (b^{2}+5a^{2})^{2}=\frac{1}{3}(80b^{2}-77a^{2}).
\end{array}\right.
\end{equation*}
Then we have
\begin{equation}\label{eq607}
b^{2}=\frac{1}{35}(108a^{4}+32a^{2}),
\end{equation}
and  $(\frac{108a^{4}+32a^{2}}{35}+5a^{2})^{2}=\frac{1}{3}(\frac{80(108a^{4}+32a^{2})}{35}-77a^{2})$.
By simplifying this equation, we have
$a^{2}(a-1)(a+1)(6a-1)(6a+1)(36a^{2}+175)=0.$
Thus $$a_{1}=a_{2}=0, a_{3}=1, a_{4}=-1, a_{5}=\frac{1}{6}, a_{6}=-\frac{1}{6}, a_{7}=\frac{5\sqrt7}{6}i, a_{8}=-\frac{5\sqrt7}{6}i.$$

Now we show that $f(1)=\pm \frac{5\sqrt7}{6}i$ is impossible. Otherwise, we have $a^{2}=-\frac{175}{36}$
and $b^{2}=\frac{2465}{36}$ by  (\ref{eq607}).

Noting that $581=24^{2}+5\times1^{2}=9^{2}+5\times10^{2}$, we have
\begin{equation}\label{eq608}f^{2}(24)=f^{2}(9)+5f^{2}(10)-5f^{2}(1)=\frac{1}{3}(575a^{4}-572a^{2}).
\end{equation}
On the other hand, we know $f(24)=f(2^{2}+5\times2^{2})=6f^{2}(2)=6b^{2}$, then by  (\ref{eq608})  we have
\begin{equation}\label{eq609}\frac{1}{3}(575a^{4}-572a^{2})-36b^{4}=0.\end{equation}

Unfortunately,  when we put $a^{2}=-\frac{175}{36}$ and $b^{2}=\frac{2465}{36}$ in the left side of (\ref{eq609}),
we obtain $\frac{1}{3}(575a^{4}-572a^{2})-36b^{4}\not=0.$
Thus it implies a contradiction and  we complete the proof.
\end{Proof}

By Lemmas \ref{lem602}-\ref{lem603} and direct calculation, we obtain the following Corollary.
\begin{cor}\label{cor604} Let $1\leq n \leq10$ and $f$ satisfy the function equation (\ref{eq601}). Then one of the following holds.

\noindent {\rm(1) } $f(n) \equiv 0$;

\noindent {\rm(2) } $f(n) =\left\{\begin{array}{ll}
 n, & \mbox { if } n \in \{6,9\}; \\
 \pm n, &  \mbox { otherwise}.
            \end{array}\right.$

\noindent {\rm(3) }$f(n) =\left\{\begin{array}{ll}
 \frac{1}{6}, & \mbox { if } n \in \{6,9\}; \\
 \pm \frac{1}{6}, &  \mbox { otherwise}.
            \end{array}\right.$
\end{cor}

Now, by Theorem \ref{them201} and Corollary \ref{cor604}, we complete the proof of Theorem \ref{them601}.

\begin{them}\label{them601}
Conjecture \ref{con202} holds for $k=5$.
\end{them}

\section{Some remarks}

\hskip.6cm So far, we know  Conjecture \ref{con202} holds for $k=1,2,3,4,5$ by \cite{2014B} and the results of Sections 3-6.
In fact, for $k=2,3,4,5$,  firstly, we find the formulae of $f(3), f(4), \ldots, f(A)$ with $f(1)$ and $f(2)$,
then evaluate the value of $f(1)$, and  show the result holds for all $n$ with $1\leq n\leq A$ by direct calculation,
finally, we show Conjecture \ref{con202} holds for $k=2,3,4,5$ by Theorem \ref{them201}.
Similarly,  we can show the cases  $k=6, 7, \ldots$ by the similar methods, so we omit it.
Of course, we expect  there exists better methods to simplify the proof.

 By the way, we find an interesting fact. The formulae (\ref{eq402})-(\ref{eq406}), (\ref{eq502}) and (\ref{eq602})-(\ref{eq603}) all hold for $k=3,4,5$.  Are they also hold for $k=6, 7$ or more? We donot know, but it is worth expecting.



\begin{thebibliography}{}
\bibitem{2013D}
A. Dubickas, P. \v{S}arka, On multiplicative functions which are additive on sums of primes, Aequat. Math., 86 (2013), 81--89.

\bibitem{2014B}
B. Ba\v{s}i\'{c}, Characterization of arithmetic functions that preserve the sum-of-squares operation, Acta Mathematica Sinica, English Series, 30 (2014), 689--695.

\bibitem{2006P}
B.M. Phong, A characterization of the identity function with an equation of Hossz\'{u} type, Publ. Math. Debrecen, 69 (2006), 219--226.

\bibitem{2004P}
B.M. Phong, On sets characterizing the identity function, Ann. Univ. Sci. Budapest. Sect. Com-put., 24 (2004), 295--306.

\bibitem{1992S}
C.A. Spiro, Additive uniqueness sets for arithmetic functions, J. Number Theory, 42 (1992), 232--246.

\bibitem{2011F}
J.H. Fang, A characterization of the identity function with equation $f(p+q+r)=f(p)+f(q)+f(r)$, Combinatorica, 31 (2011), 697--701.

\bibitem{1997D}
J.-M. De Koninck, K\'{a}tai I. and B.M. Phong, A new characteristic of the identity function, J. of Number Theory, 63 (1997), 325--338.

\bibitem{2006I}
K.-H. Indlekofer, B. Phong, Additive uniqueness sets for multiplicative functions, Ann. Univ. Sci. Budapest. Sect. Comput., 26 (2006), 65--77.

\bibitem{1996C}
P.V. Chung, Multiplicative functions satisfying the equation $f(u^{2}+v^2)=f^{2}(u)+f^{2}(v)$, Math. Slovaca, 46 (1996), 165--171.

\bibitem{2016C}
Y.G. Chen, J.H. Fan, P.Z. Yuan, Y.P. Zheng, On multiplicative functions with $f(p+q+n_{0})=f(p)+f(q)+f(n_{0})$, J. Number Theory, 165 (2016), 01--20.


\end{thebibliography}
\end{document}